\ifx\crversion\undefined
\documentclass[11pt]{amsart}
\else 
\documentclass[CRMATH,Unicode,manuscript]{cedram}
\fi 

\usepackage{amsmath,amsthm,amssymb,amsfonts}
\usepackage{xcolor}

\usepackage{xargs}  

\usepackage{hyperref}
\usepackage{enumitem}
\setenumerate{leftmargin=*}
\setitemize{leftmargin=*}
\usepackage{amscd}

\numberwithin{equation}{section}

\swapnumbers

\ifx\crversion\undefined
\newtheorem{thm}[subsection]{Theorem}

\newtheorem{proposition}[subsection]{Proposition}

\newtheorem{conjecture}{Conjecture}
\fi

\newcounter{consta}

\newcounter{constk}
\renewcommand{\theconstk}{{\kappa_{\arabic{constk}}}}
\newcommand{\constk}{\refstepcounter{constk}\theconstk}

\newcounter{constc}

\newcounter{constE}
\renewcommand{\theconstE}{{{C}_{\arabic{constE}}}}
\newcommand{\constE}{\refstepcounter{constE}\theconstE}

\def\bbr{\mathbb{R}}

\def\bbc{\mathbb{C}}

\def\bbn{\mathbb{N}}
\def\Gbf{\mathbb{G}}

\def\R{\bbr}
\def\C{\bbc}

\def\hfrak{\mathfrak{h}}

\def\rfrak{\mathfrak{r}}

\def\gfrak{\mathfrak{g}}

\def\Gbf{\mathbf{G}}

\def\G{\Gbf}

\DeclareMathOperator\diff{d}

\def\vol{{\rm{vol}}}
\def\SL{{\rm{SL}}}

\def\SO{{\rm{SO}}}
\def\Lie{{\rm Lie}}

\DeclareMathOperator\Ad{Ad}

\def\vare{\varepsilon}

\def\zg0{Z_{G_\omega}(s)}

\def\zg{Z_G(s)}

\def\be{\begin{equation}}
\def\ee{\end{equation}}

\def\dist{{\rm dist}}

\def\lf{\mathfrak l}
\def\Sob{{\mathcal S}}

\def\diff{\operatorname{d}}

\def\dist{d}

\def\mixexp{\kappa_0}

\def\boxH{\mathsf B^H}
\def\boxHs{\mathsf B^{H}}
\def\boxG{\mathsf B^{G}}

\def\rwm{\nu}

\def\rws{t}

\def\dexp{\delta}
\def\rel{r}

\def\ave{\int_{0}^1}
\def\uvk{u_\rel}
\def\uvkd{\diff\!\rel}

\def\mfm{\hat f}

\def\margI{I}
\def\noI{\psi}

\def\lf{\mathbb R}
\def\qlf{\mathbb C}

\def\qi{i}

\def\nuni{e}
\def\coneH{\mathsf E}
\def\cone{\mathcal E}

\newcommand{\scmf}{\varrho}
\newcommand{\sfh}{\mathsf h}

\makeatletter
\newcommand*\bigcdot{\mathpalette\bigcdot@{.5}}
\newcommand*\bigcdot@[2]{\mathbin{\vcenter{\hbox{\scalebox{#2}{$\m@th#1\bullet$}}}}}
\makeatother

\def\XXint#1#2#3{{\setbox0=\hbox{$#1{#2#3}{\int}$ }
\vcenter{\hbox{$#2#3$ }}\kern-.6\wd0}}

\newcommand{\balpha}{\alpha} 
\newcommand{\trct}{\mathsf K}

\begin{document}
\title[Polynomial effective equidistribution]{Polynomial effective equidistribution}

\ifx\crversion\undefined
\author{E.~Lindenstrauss}
\address{E.L.: The Einstein Institute of Mathematics, Edmond J.\ Safra Campus, 
Givat Ram, The Hebrew University of Jerusalem, Jerusalem, 91904, Israel}
\email{elon@math.huji.ac.il}
\thanks{E.L.\ acknowledges support by ERC 2020 grant HomDyn (grant no.~833423).}

\author{A.~Mohammadi}
\address{A.M.: Department of Mathematics, University of California, San Diego, CA 92093}
\email{ammohammadi@ucsd.edu}
\thanks{A.M.\ acknowledges support by the NSF, grants DMS-1764246 and 2055122.}

\author{Z.~Wang}
\address{Z.W.: Pennsylvania State University,
Department of Mathematics, University Park, PA 16802}
\email{zhirenw@psu.edu}
\thanks{Z.W.\ acknowledges support by the NSF, grant  DMS-1753042.}
\else
\author{\firstname{Elon} \lastname{Lindenstrauss}}
\address{E.L.: The Einstein Institute of Mathematics, Edmond J. Safra Campus, 
Givat Ram, The Hebrew University of Jerusalem, Jerusalem, 91904, Israel}
\email{elon@math.huji.ac.il}
\thanks{E.L.\ acknowledges support by ERC 2020 grant HomDyn (grant no.~833423).}

\author{\firstname{Amir} \lastname{Mohammadi}}
\address{A.M.: Department of Mathematics, University of California, San Diego, CA 92093}
\email{ammohammadi@ucsd.edu}
\thanks{A.M.\ acknowledges support by the NSF, grants DMS-1764246 and 2055122.}

\author{\firstname{Zhiren} \lastname{Wang}}
\address{Z.W.: Pennsylvania State University,
Department of Mathematics, University Park, PA 16802}
\email{zhirenw@psu.edu}
\thanks{Z.W.\ acknowledges support by the NSF, grant  DMS-1753042.}

\fi

\begin{abstract}
We prove effective equidistribution theorems, with polynomial error rate, for orbits of the unipotent subgroups of $\SL_2(\lf)$ 
in arithmetic quotients of $\SL_2(\qlf)$ and $\SL_2(\lf)\times\SL_2(\lf)$. 

The proof is based on the use of a Margulis function, tools from incidence geometry, 
and the spectral gap of the ambient space. 
\end{abstract}

\ifx\crversion\undefined
\else
\begin{altabstract} 
Nous prouvons des th\'{e}or\'{e}mes d'\'{e}quidistribution effectifs, avec un taux d'erreur polynomial pour les orbites des sous-groupes unipotents de $\SL_2(\lf)$ en quotients arithm\'{e}tiques de $\SL_2(\qlf)$ et $\SL_2(\lf)\times\SL_2(\lf)$.

La preuve est bas\'{e}e sur l'utilisation d'une fonction de Margulis, des outils de la g\'{e}om\'{e}trie d'incidence,
et le trou spectral de l'espace ambiant.
\end{altabstract}
\fi
\maketitle

\section{Introduction}
A landmark result of Ratner~\cite{Ratner-topological} states that if $G$ is a Lie group, $\Gamma$ a lattice in $G$ and if $u_t$ is a one-parameter $\Ad$-unipotent subgroup of $G$, then for {\bf any} $x \in G/\Gamma$ the orbit $u_t.x$ is equidistributed in a periodic orbit of some subgroup $L <G$ that contains both the one parameter group $u_t$ and the initial point $x$. We say an orbit $L.x$ of a  group $L$ in some space $X$ is periodic if the stabilizer of $x$ in $L$ is a lattice in $L$, equivalently that the stabilizer of $x$ in $L$ is discrete and $L.x$ supports a unique $L$-invariant probability measure $m_{L.x}$;
and $u_t.x$ is equidistributed in $L.x$ in the sense that 
\begin{equation}\label{eq:equidist}
\frac1T \int_{0}^T f(u_t.x) dt \to \int f dm_{L.x} \qquad\text{for any $f \in C_0(G/\Gamma)$.}
\end{equation}
In order to prove this equidistribution result, Ratner first classified the $u_t$-invariant probability measures on $G/\Gamma$ \cite{Ratner-Acta, Ratner-measure}; the proof also uses the non-divergence properties of unipotent flows established by Dani and Margulis \cite{Marg-Nondiv, Dani-Nondiv-1, Dani-Nondiv-2}.

A key motivation behind Ratner's equidistribution theorem for one parameter unipotent flows has been to establish Raghunathan's conjecture regarding the possible orbit closures of groups generated by one-parameter unipotent groups; using the equidistribution theorem Ratner proved that if $G$ and $\Gamma$ are as above, and if $H<G$ is generated by one parameter unipotent groups, then for any $x \in G/\Gamma$ one has that $\overline{H.x}=L.x$ where $H < L < G$ and $L.x$ is periodic. Important special cases of Raghunathan's conjecture were proven earlier by Margulis and by Dani and Margulis using a different more direct approach, which in particular gave a proof of a rather strong form of the longstanding Oppenheim conjecture \cite{Margulis-Oppenheim,DM-Oppenheim,DM-MathAnn}.

These results have had a surprisingly rich trove of applications in number theory and beyond. One drawback of this method, compared to more traditional number theoretic techniques, is that these equidistribution results were neither effective nor quantitative. Indeed, Ratner's proof relies heavily on the pointwise ergodic theorem and Lusin's theorem, both rather fundamental theorems but that do not have good effective analogues. The Dani-Margulis method is somewhat easier to make effective and a result in this direction was given by Margulis and the first named author in \cite{LM-Oppenheim}; moreover a general result in this direction was announced by Margulis, Shah and two of us (E.L. and A.M.) with the first installment of this work appearing in \cite{LMMS}. However, this only gives density properties of these flows, not equidistribution, and the rate of density obtained in this way is far from optimal (polylog at best, though in most cases one only has an iterative-log type bound).

\medskip

In this paper we announce a quantitative equidistribution result for orbits of a one parameter unipotent group on quotients $G/\Gamma$ where $G$ is either $\SL_2(\qlf)$ or $\SL_2(\lf)\times\SL_2(\lf)$ with a polynomial error rate, which is the first quantitative equidistribution statement for individual orbits of unipotent flows on quotients of semi-simple groups beyond the horospherical case. Our approach builds on the paper \cite{LM-PolyDensity} by the first two authors, where an effective density result with a polynomial rate for orbits of a Borel subgroup of a subgroup $H \simeq \SL_2(\bbr)$ of $G$ was proved.

\medskip

Recall that a group $N<G$ is horospheric if there is some $g \in G$ so that 
\[
N = \{h \in G: g^{-n} h g^{n} \to 1 \text{ as $n\to\infty$}\}.
\]
For instance, the one parameter unipotent group
\[
\left\{\begin{pmatrix}1&r\\0&1\end{pmatrix}: r\in\bbr\right\}
\]
is horospheric in $\SL_2(\R)$ as are the groups
\[
\left\{\begin{pmatrix}1&r+is\\0&1\end{pmatrix}:r,s\in\bbr \right \}
\quad\text{and}\quad
\left\{\left(\begin{pmatrix}1&r\\0&1\end{pmatrix}, \begin{pmatrix}1&s\\0&1\end{pmatrix}\right): r,s\in\bbr\right\}
\]
in $\SL_2(\C)$ and $\SL_2(\R)\times\SL_2(\R)$, respectively.
The classification of invariant measures and orbit closures for horospherical flows was established prior to Ratner's work by Hedlund, Furstenberg, Dani, Veech and others, and this has been understood for some time also quantitatively since one can relate the distribution properties of individual $N$ orbits to the ergodic theoretic properties of the action of $g$ on $G/\Gamma$ (cf.~\S\ref{large dim=>equidistribution} for more details).

\medskip

We also mention that while our result is the first quantitative equidistribution statement for individual orbits of unipotent flows on quotients of semi-simple groups beyond the horospherical case\footnote{As pointed out by Venkatesh in \cite{Venkatesh-Sparse}, the distinction between horospheric and non-horospheric is not completely clear cut, and indeed the results of that paper can also be recast as a nonhorospheric equidistribution problem; cf.\ \cite{FFT-twisted} as well as \cite{Katz-Quantitative}. There are also some quantitative equidistribution results for particular types of unipotent orbits, e.g.\ \cite{Chow-Yang} by Chow and Lei Yang.} there have been some other interesting quantitative equidistribution results. When $G$ is unipotent, an effective equidistribution result for unipotent flows on $G/\Gamma$ was given by Green and Tao in \cite{GrTao-Nil}, and was a key ingredient in a series of works by Green, Tao and Ziegler about linear equations in primes. In the case of quotients of the skew product $G=\SL_2(\R)\ltimes \R^2$, Strombergsson \cite{Strom-Semi} has an effective equidistribution result for one parameter unipotent orbits (which are not horospheric in $G$, but project to a horospheric group on $\SL_2(\R)$), and this has been generalized by several authors, in particular by Wooyeon Kim \cite{ASLn-Kim} (using a completely different argument) to $\SL_n(\R)\ltimes\R^n$. 
Moreover there is an important work of Einsiedler, Margulis and Venkatesh \cite{EMV} showing that periodic orbits of semisimple subgroups $H$ of a semisimple group $G$ are quantitatively equidistributed in an appropriate homogeneous subspace of $G/\Gamma$ if $\Gamma$ is a congruence lattice and $H$ has finite centralizer in $G$. Subsequently Einsiedler, Margulis, Venkatesh and the second named author by using Prasad's volume formula and a more adelic view point were able to prove such an equidistribution result for periodic orbits of maximal semisimple subgroups of $G$ when the subgroup is allowed to vary \cite{EMMV} with arithmetic applications (the equidistribution of periodic orbits of semisimple groups is also closely connected to the equidistribution of Hecke points; a quantitative treatment of such equidistribution was given by Clozel, Oh and Ullmo in \cite{Clozel-Oh-Ullmo}).

In a slightly different direction, Bourgain, Furman, Mozes and the first named author \cite{BFLM} gave a quantitative equidistribution result with exponential rates (this is analogous to polynomial rates in our problem) for random walks by automorphisms of the torus. In this case this equidistribution result was new even without rates. There have been several extensions of this result, in particular \cite{He-Saxce} where a proximality assumption was removed.  Kim in \cite{ASLn-Kim} used the techniques of \cite{BFLM} to study  $\SL_n(\R)\ltimes\R^n$. Our work is also heavily influenced by \cite{BFLM}.

\bigskip

We now proceed to describe our results and some of the ingredients of the proofs.
Let     
\[
G=\SL_2(\qlf)\quad\text{ or }\quad G=\SL_2(\lf)\times\SL_2(\lf).
\] 
Let $\Gamma\subset G$ be a lattice, and put $X=G/\Gamma$. We let $m_X$ denote the  $G$-invariant probability measure on $X$.   
Throughout the paper, we will denote by $H$ a subgroup of $G$ isomorphic to $\SL_2(\lf)$, namely
\[
\SL_2(\lf)\subset\SL_2(\qlf)\quad\text{or}\quad\{(g,g): g\in\SL_2(\lf)\}\subset \SL_2(\lf)\times\SL_2(\lf).
\]

For all $t, r\in\bbr$, let $a_t$ and $u_r$ denote the image of 
\[
\begin{pmatrix} e^{t/2} & o \\ 0 & e^{-t/2}\end{pmatrix}\quad\text{and}\quad \begin{pmatrix} 1 & r \\ 0 & 1\end{pmatrix},
\]
in $H$, respectively. 

We fix maximal compact subgroups ${\rm SU}(2)\subset \SL_2(\bbc)$ and $\SO(2)\times\SO(2)\subset\SL_2(\bbr)\times\SL_2(\bbr)$.
Let $\dist$ be the right invariant metric on $G$ which is defined using the Killing form and the aforementioned maximal compact subgroups. This metric induces a metric $\dist_X$ on $X$, 
and natural volume forms on $X$ and its submanifolds. We define the injectivity radius of a point $x\in X$ using this metric. 
In the sequel, $\|\;\|$ denotes the maximum norm on ${\rm Mat}_2(\qlf)$ or 
${\rm Mat}_{2}(\lf)\times {\rm Mat}_{2}(\lf)$ with respect to the standard basis. 

\medskip

The following are the main results in this paper.

\begin{thm}\label{thm: main}\label{thm:main}
Assume $\Gamma$ is an arithmetic lattice. 
For every $x_0\in X$ and large enough $R$ (depending explicitly on the injectivity radius of $x_0$), for any $T \geq R$, at least one of the following holds.
\begin{enumerate}
\item For every $\varphi\in C_c^\infty(X)$, we have 
\[
\biggl|\int_0^1 \varphi(a_{\log T}u_rx_0)\diff\!r-\int \varphi\diff\!m_X\biggr|\leq \ref{c:main-1}\Sob(\varphi)R^{-\ref{k:main-1}}
\]
where $\Sob(\varphi)$ is a certain Sobolev norm. 
\item There exists $x\in X$ such that $Hx$ is periodic with $\vol(Hx)\leq R$, and 
\[
\dist_X(x,x_0)\leq \ref{c:main-1}R^{A}T^{-1}.
\] 
\end{enumerate} 
The constants $A$, $\constk\label{k:main-1}$, and $\constE\label{c:main-1}$ are positive, and depend on $X$ but not on $x_0$.
\end{thm}


Theorem~\ref{thm:main} can be viewed as an effective version of \cite[Thm.~1.4]{Shah-Expanding}. Combining Theorem~\ref{thm:main} and the Dani--Margulis linearization method \cite{DM-Linearization} (cf.\ also Shah \cite{Shah-MathAnn}),
that allows to control the amount of time a unipotent trajectory spends near
invariant subvarieties of a homogeneous space, 
we also obtain an effective equidistribution theorem for long pieces of unipotent orbits (more precisely, we use a sharp form of the linearization method taken from~\cite{LMMS}).

\begin{thm}\label{thm:main unipotent}
Assume $\Gamma$ is an arithmetic lattice.
For every $x_0\in X$ and large enough $R$ (depending explicitly on the injectivity radius of $x_0$), for any $T\geq R$, at least one of the following holds.
\begin{enumerate}
\item For every $\varphi\in C_c^\infty(X)$, we have 
\[
\biggl|\frac{1}{T}\int_0^{T} \varphi(u_rx_0)\diff\!r-\int \varphi\diff\!m_X\biggr|\leq \ref{c:main uni}\Sob(\varphi)R^{-\ref{k:main uni}}
\]
where $\Sob(\varphi)$ is a certain Sobolev norm. 
\item There exists $x \in G/\Gamma$ with $\vol(H.x)\leq R$, 
and for every $r\in [0,T]$ there exists $g\in G$ with $\|g\|\leq R$ so that  
\[
\dist_X(u_{s}x_0, gH.x)\leq \ref{c:main uni}R^{A_1}\left(\frac{|s-r|} T\right)^{1/A_2}\quad\text{for all $s\in[0,T]$.}
\] 
\item For every $r\in[0,T]$ and $t \in [\log R, \log T]$, the injectivity radius at $a_{-t}u_r x_0$ is at most $\ref{c:main uni} R^{A_1}e^{-t/A_2}$.

\end{enumerate} 
The constants $A_1$, $A_2$, $\constk\label{k:main uni}$, and $\constE\label{c:main uni}$ are positive, and depend on $X$ but not on $x_0$. 
\end{thm}

The assumption in Theorem~\ref{thm:main} that $\Gamma$ is arithmetic may be relaxed. 
Let us say $\Gamma$ has {\em algebraic entries} if the following is satisfied:
there is a number field $F$, a semisimple $F$-group $\G$ of adjoint type, 
and a place $v$ of $F$ so that $F_v=\bbr$ and $\G(F_v)$ and $G$ are locally isomorphic --- in which case there is a surjective homomorphism from $G$ onto the connected component of the identity in $\G(F_v)$ ---
and the image of $\Gamma$ in $\G(F_v)$ (possibly after conjugation) is contained in $\G(F)$. Every arithmetic lattice has algebraic entries, but there are lattices with algebraic entries that are not arithmetic.

Note that the condition that  $\Gamma$ has {\em algebraic entries} is automatically satisfied if $\Gamma$ is an irreducible lattice in $\SL_2(\bbr)\times \SL_2(\bbr)$ or if $G=\SL(2,\C)$. Indeed, by arithmeticity theorems of Selberg and Margulis, irreducible lattices in $\SL_2(\bbr)\times \SL_2(\bbr)$ are arithmetic~\cite[Ch.~IX]{Margulis-Book}. Moreover, by local rigidity, lattices in $\SL_2(\bbc)$ always have algebraic entries \cite[Thm.~0.11]{GarRag-FD} (see also~\cite{Selb, Weil-1, Weil-2}).

\begin{thm}\label{thm: main arithmeticity relaxed}
Assume $\Gamma$ is a lattice which has algebraic entries.
For every $0<\dexp<1$, every $x_0\in X$ and large enough $T$ (depending explicitly on $\dexp$ and the injectivity radius of $x_0$) 
at least one of the following holds.
\begin{enumerate}
\item For every $\varphi\in C_c^\infty(X)$, we have 
\[
\biggl|\int_0^1 \varphi(a_{\log T}u_rx_0)\diff\!r-\int \varphi\diff\!m_X\biggr|\leq \ref{c:main-1}\Sob(\varphi)T^{-\dexp^2\ref{k:main-1}}
\]
where $\Sob(\varphi)$ is a certain Sobolev norm. 
\item There exists $x\in X$ with
\[
\dist_X(x,x_0)\leq \ref{c:main-1}T^{-1/A},
\]
 satisfying the following: there are elements $\gamma_1$ and $\gamma_2$ in ${\rm Stab}_H(x)$
with $\|\gamma_i\|\leq T^{\dexp}$ for $i=1,2$ so that the group generated by $\{\gamma_1,\gamma_2\}$ is Zariski dense in $H$.
\end{enumerate} 
The constants $A$, $\ref{k:main-1}$, and $\ref{c:main-1}$ are positive, and depend on $X$ but not on $\delta$ and $x_0$.
\end{thm}

The obstacle to effective equidistribution in Theorem~\ref{thm:main} is much cleaner and simpler than in Theorem~\ref{thm:main unipotent}. This is not an artifact of the proof but a reflection of reality; a unipotent orbit may fail to equidistribute at the expected rate without it staying near a single period orbit of some subgroup $\{u_t\}< L< G$: one must allow a slow drift of the periodic orbit in the direction of the centralizer of $u_t$. Unlike the work of Shah in \cite{Shah-Expanding}, where (in particular) a non-effective version of Theorem~\ref{thm:main} is proved relying on Ratner's measure classification theorem for unipotent flows, our proof goes the other way, first establishing Theorem~\ref{thm:main}, and then reducing Theorem~\ref{thm:main unipotent} from it using a linearization and non-divergence argument. 

An extremely interesting analogue to unipotent flows on homogeneous spaces is given by the action of $\SL_2(\R)$ and its subgroups on strata of abelian differentials. 
Let $g \geq 1$, and 
let $\alpha = (\alpha_1,\dots, \alpha_n)$ be a partition of $2g-2$. Let $\mathcal H(\alpha)$ be the corresponding stratum of abelian differentials, i.e., the space of pairs $(M,\omega)$ where $M$ is a compact Riemann surface with genus $g$ and $\omega$ is a holomorphic $1$-form on $M$ whose zeroes have
multiplicities $\alpha_1, \dots, \alpha_n$. The form $\omega$ defines a canonical flat metric on $M$ with conical singularities and a natural area form. Let $\mathcal H_1(\alpha)$ be the space of unit area surfaces in $\mathcal H(\alpha)$. The space $\mathcal H(\alpha)$
admits a natural action of $\SL_2(\bbr)$; this action preserves the unit area hyperboloid $\mathcal H_1(\alpha)$.

A celebrated theorem of Eskin and Mirzakhani~\cite{EMir-Meas} shows that any $P$-invariant ergodic measure is $\SL_2(\bbr)$-invariant and is supported on an affine invariant manifold, where $P$ denotes the group of upper triangular matrices in $\SL_2(\bbr)$. We shall refer to these measures as {\em affine invariant measures}. 
Moreover if we define, for any interval $I \subset \R$ and $x \in \mathcal H_1(\alpha)$, the probability measure $\mu_I^x$ on $\mathcal H_1(\alpha)$ by
\[
\mu_{I}^{x}=|I|^{-1}\int_I \delta_{u_s x}\diff\! s,
\]
then Eskin, Mirzakhani and the second named author \cite{EMM-Orbit} showed that for any  $x \in \mathcal H_1(\alpha)$ the limit
\begin{equation}\label{eq:EMM-Orbit}
    \lim_{T\to\infty} \frac{1}{T}\int_{t=0}^T a_t\mu_{[0,1]}^x \diff\! t \qquad\text{exists in weak${}^*$ sense}
\end{equation}
and is equal to an ($\SL(2,\R)$-invariant) affine invariant probability measure with $x$ in its support.
On the other hand, there are several results, in particular by Chaika, Smillie and B.~Weiss in~\cite{CSW-Tremors}, that show that an analogue of Ratner's equidistribution theorem (or our Theorem~\ref{thm:main unipotent}) fails to hold in this setting, for instance for some $x$ the sequence of measure $\mu_{[0,T]}^{x}$ may fail to converge as $T\to\infty$, or may converge to a non-ergodic measure. We believe the following strengthening of \eqref{eq:EMM-Orbit}, which we have learned has already been conjectured by Forni in \cite[Conj.~1.4]{Forni-limits}, should however hold:

\begin{conjecture}
Let $\mathcal H_1(\alpha)$ be the space of unit area surfaces in stratum of abelian differentials on a genus $g$ surface whose zeros have multiplicities given by $\alpha=(\alpha_1, \dots, \alpha_n)$, and let $x\in\mathcal H_1(\alpha)$.
Then $\lim_{t\to\infty} a_t\mu^{x}_{[0,1]}$ exists in the weak${}^*$ sense and is equal to an affine invariant measure with $x$ in its support.  
\end{conjecture}

\noindent Of course, once one establishes that $\lim_{t\to\infty} a_t\mu^{x}_{[0,1]}$ exists, the rest follows from \cite{EMM-Orbit}.

\section{Some preliminaries}
We discuss the proofs of Theorems~\ref{thm:main} and~\ref{thm: main arithmeticity relaxed}. As mentioned above, Theorem~\ref{thm:main unipotent} is proved by combining Theorem~\ref{thm:main} and the linearization techniques;  in this announcement we focus on the proof of the former results. We note that the idea of using equidistribution of expanding translates of a fixed piece of a $U$ orbit of the type $\{a_t u_s .x: 0\leq s \leq 1\}$ to deduce equidistribution of a large segment of a non-translated $U$ orbit $\{u_s.x: 0\leq s \leq T\}$ is quite classical. 

\medskip

Let $U\subset N$ denote the group of upper triangular unipotent matrices in $H\subset G$, respectively, 
and let $A=\{a_t: t\in\bbr\}\subset H$.
More explicitly, if $G=\SL_2(\qlf)$, then  
\[
N=\left\{n(\rel,s)=\begin{pmatrix} 1 & \rel+is\\ 0 &1\end{pmatrix}: (r,s)\in\lf^2\right\}
\] 
and $U=\{n(\rel,0):\rel\in \lf\}$; note that $n(r,0)=u_\rel$ for $\rel\in\bbr$. Let 
\[
V=\{n(0, s)=v_s: s\in \lf\};
\]
if $G=\SL_2(\lf)\times \SL_2(\lf)$, then  
\[
N=\left\{n(\rel,s)=\left(\begin{pmatrix} 1 & \rel+s\\ 0 &1\end{pmatrix},\begin{pmatrix} 1 & \rel\\ 0 &1\end{pmatrix}\right): (r,s)\in\lf^2\right\}
\] 
and $U=\{n(\rel,0):\rel\in \lf\}$.
As before, $n(\rel,0)=u_r$ for $\rel\in\bbr$. Let \[V=\{n(0,s)=v_s: s\in\lf\}.\] In both cases, we have $N=UV$.
Let us denote the transpose of $U$ by $U^-$ and its elements by $u^-_r$.

\medskip

Let $\gfrak=\Lie(G)$, that is, $\gfrak=\mathfrak {sl}_2(\qlf)$ or $\gfrak=\mathfrak{sl}_2(\lf)\oplus \mathfrak{sl}_2(\lf)$.
Let $\rfrak=\qi\mathfrak{sl}_2(\bbr)$ if $\gfrak=\mathfrak{sl}_2(\qlf)$ and 
$\rfrak=\mathfrak{sl}_2(\bbr)\oplus\{0\}$
if $\gfrak=\mathfrak{sl}_2(\bbr)\oplus \mathfrak{sl}_2(\bbr)$.
In either case $\mathfrak g=\mathfrak h\oplus \mathfrak r$ where $\hfrak=\Lie(H)\simeq\mathfrak{sl}_2(\bbr)$, and both $\hfrak$ and $\rfrak$ are $\Ad(H)$-invariant.

We fix a norm on $\hfrak$ by taking the maximum norm where the coordinates are given by $\Lie(U)$, $\Lie(U^-)$ and $\Lie(A)$; similarly fix a norm on $\rfrak$. 
By taking maximum of these two norms, we obtain a norm on $\gfrak$. All these norms will be denoted by $\|\;\|$.

For all ${\beta}>0$, we define $\boxG_{\beta}:=\exp(B_\hfrak(0,\beta))\cdot\exp(B_\rfrak(0,\beta))$ 
where $B_\bullet(0,{\beta})$ denotes the ball of radius $\beta$ in $\bullet$ with respect to $\|\;\|$.

We also define
\[
\boxH_{\beta}:=\{u_s^-:|s|\leq {\beta}\}\cdot\{a_t: |t|\leq \beta\}\cdot\{u_\rel:|\rel|\leq {\beta}\}
\] 
for all $0<{\beta}<1$.

\medskip

For the sake of simplicity of the exposition here, let us assume $X$ is compact. 
Let $0<\eta_0<1$ be so that the map 
$g\mapsto gx$ is injective on $\boxG_{100\eta_0}$ for all $x\in X$.

\section{From large dimension to equidistribution}\label{large dim=>equidistribution}
We begin with an equidistribution theorem which is of independent interest.  
In the proof of Theorem~\ref{thm:main}, this proposition will be applied to each of the sets obtained in the bootstrap phase, see Proposition~\ref{propos: main bootstrap equi}. 

\medskip

Let us recall the following quantitative decay of correlations for the ambient space $X$: 
There exists $0<\mixexp\leq1$ so that  
\be\label{eq:actual-mixing-intro}
\biggl|\int \varphi(gx)\psi(x)\diff\!{m_X}-\int\varphi\diff\!{m_X}\int\psi\diff\!{m_X}\biggr|\ll \Sob(\varphi)\Sob(\psi) \nuni^{-\mixexp d(e,g)}
\ee
for all $\varphi,\psi\in C^\infty_c(X)+\bbc\cdot 1$, where $m_X$ is the $G$-invariant probability measure on $X$ and $d$ is our fixed right $G$-invariant metric on $G$. See, e.g., \cite[\S2.4]{KMnonquasi} and references there for~\eqref{eq:actual-mixing-intro}; we note that $\mixexp$ is absolute if $\Gamma$ is a congruence subgroup. This is known in much greater generality, but the cases relevant to our paper are due to Selberg and Jacquet-Langlands~\cite{Selberg-ThreeSixteenth, Jacquet-Langlands}.

The quantitative decay of correlation can be used to establish quantitative results regarding the equidistribution of translates of pieces of an $N$-orbit. Specifically we employ the results in~\cite{KMnonquasi}, but there is rich literature around the subject; a more complete list of such papers can be found in \cite[\S1]{LM-PolyDensity}.

Now let $\xi:[0,1]\to \rfrak$ be a smooth nonconstant curve. Then using the quantitative results regarding equidistribution of translates of pieces of an $N$-orbit such as~\cite{KMnonquasi}, one can show that for every $x\in X$,
\[
a_\tau\Bigl\{u_r\exp(\xi(s)).x:r,s\in [0,1]\Bigr\}
\]
is equidistributed in $X$ as $\tau\to\infty$ (with a rate which is polynomial in $\nuni^{-\tau}$).  
The key point in the deduction of this equidistribution result from the equidistribution of shifted $N$ orbits is that conjugation by $a_\tau$ moves $u_r\exp(\xi(s))$ to the direction of 
$N$, hence the above average essentially reduces to an average on a $N$ orbit.

\medskip 

Roughly speaking, the following proposition states that one may replace the curve $\{\xi(s):s\in[0,1]\}$ with a measure on 
$\rfrak$ so long as the measure has dimension $\geq 1-\theta$, for an appropriate choice of $\theta$ depending on $\mixexp$. 

\medskip

The precise formulation is the following.

\begin{proposition}\label{prop: equi translates of cone intro}
For any $\theta>0$ and $c>0$ there is a $\constk\label{k: rho/varrho exp}$ so that the following holds:
Let $0<\varrho<10^{-6}$, and let $F\subset B_\rfrak\Bigl(0,\varrho\Bigr)$ 
be a finite set satisfying  
\[
\frac{\#(F\cap B_\rfrak(0,b))}{\#F}\leq \rho^{-c}\Bigl(b/\varrho\Bigr)^{1-\theta}\quad\text{for all $b\geq \rho$}
\]
where $\rho <\varrho^{20}$.

Then for all $x\in X$ and all $2\log(1/\varrho)\leq \tau\leq \frac{1}{10} \log(1/\rho)$, we have 
\begin{equation*}
   \biggl|\int_0^1\frac{1}{\#F}\sum_{w\in F}\varphi(a_\tau u_r\exp(w)x)\diff\!r-\int\varphi\diff\!m_X\biggr|\ll_X 
   \Sob(\varphi)\max\left((\rho/\varrho)^{\ref{k: rho/varrho exp}},\rho^{-2c} e^{2\tau\theta}\varrho^{\kappa_0^2/M}\right), 
\end{equation*}
where $\Sob(\varphi)$ is a certain Sobolev norm and $M\geq 100$ an absolute constant. 
\end{proposition}

The proof of this proposition is significantly more delicate than that of the ``toy version" of a shifted curve, and relies on an adaptation of a projection theorem due to K\"{a}enm\"{a}ki, Orponen, and Venieri~\cite{kenmki2017marstrandtype}, based on the works of 
Wolff~\cite{Wolff}, Schlag~\cite{Schlag}, and~\cite{Zahl}, in conjunction with a sparse equidistribution argument due to Venkatesh~\cite{Venkatesh-Sparse}. 
These elements also played a crucial role in previous work by E.L.\ and A.M.~\cite{LM-PolyDensity} regarding quantitative density for the action of $AU$ on the spaces we consider here.

\medskip
The goal in the remaining steps is to show that unless Theorem~\ref{thm:main}(2) holds, we can find a subset $J\subset[0,1]$, where $[0,1]\setminus J$ has an exponentially small measure, and up to an exponentially small error the uniform measure on $\{a_{A\log T}u_r x_0: r\in J\}$
can be decomposed into counting measures on sets each of which satisfies the conditions in Proposition~\ref{prop: equi translates of cone intro}.

\section{Inheritance of the Diophantine property}
If part~(2) in Theorem~\ref{thm:main} holds, we are done. 
Let us, therefore, assume the alternative, which gives some Diophantine condition on the point $x_0$ in terms of its distance to nearby periodic $H$-orbits.
The first step in the proof is to improve this Diophantine condition, perhaps not at $x_0$, but at some (indeed all except a set of exponentially small measure) point on the translation of a the $U$-orbit segment $\{u_r.x_0: r \in [0,1]\}$ by a big diagonal element $a_{t}$.

\begin{proposition}\label{prop: linearization translates}
There exist $D_0$ (absolute) and $\constE\label{c: linear trans}$ (depending on $X$) so that the following holds. 
Let $S\geq R$, and let $x_0\in X$ be so that 
\[
\dist_X(x,x_0)\geq S^{-1}
\] 
for all $x$ with $\vol(Hx)\leq R$. Then for all $s\geq \log S$ we have the following 
\[
\biggl|\biggl\{r\in [0,1]: \begin{array}{c}\text{There is $x$ with } \vol(Hx)\leq R\\ 
\text{so that }\dist_X(x, a_{s}u_rx_0)\leq \ref{c: linear trans}R^{-D_0}\end{array}\biggr\}\biggr|\leq \ref{c: linear trans}R^{-1}.
\]
\end{proposition}


In the proof of Proposition~\ref{prop: linearization translates} we consider each periodic orbit individually, and then use the fact that the number of periodic $H$-orbits with volume $\leq R$ 
in $X$ is $\ll R^6$, see e.g.~\cite[\S10]{MO-MargFun} to conclude.
The desired result for an individual orbit can be proved using Margulis functions for periodic $H$-orbits similar to those which were used in~\cite[\S9]{LM-PolyDensity}, see also~\cite[Prop.\ 2.13]{EMM-Orbit}.

It is worth mentioning that even though~\cite[Thm.~1.4]{LMMS} concerns long pieces of $U$-orbits and Proposition~\ref{prop: linearization translates} concerns translates of pieces of $U$-orbits, similar tools are applicable here as well. In particular, a version of Proposition~\ref{prop: linearization translates} can also be proved using the methods of~\cite{LMMS}.

\section{Closing lemma}
Let $t>0$ be a large parameter and fix some $\nuni^{-0.01t}<\beta=\nuni^{-\kappa t}<\eta_0$; in our application, $\kappa$ will be chosen to be $\ll 1/D_0$ where the implied constant depends on $X$ and $D_0$ is as in Proposition~\ref{prop: linearization translates}. 

For every $\tau\geq0$, put
\[
\coneH_{\tau}=\boxHs_{\beta}\cdot a_\tau\cdot \{u_r: r\in [0,1]\} \subset H.
\]

If $y\in X$ is so that the map $h\mapsto hy$ is injective on $\coneH_{\tau}$, then $\mu_{\coneH_\tau.y}$ 
denotes the pushforward of the normalized Haar measure on $\coneH_\tau$ to $\coneH_\tau.y\subset X$.

Let $\tau\geq 0$ and $y\in X$. For every $z\in\coneH_\tau.y$, put
\[
\margI_\tau(z):=\Bigl\{w\in \rfrak: 0<\|w\|<\eta_0 \text{ and } \exp(w) z\in \coneH_\tau.y\Bigr\};
\]
this is a finite subset of $\rfrak$ since $\coneH_\tau$ is bounded ---  
we will define $\margI_\cone(h,z)$ 
for all $h\in H$ and more general sets $\cone$ in the bootstrap phase below.

Let $0<\alpha<1$. Define the function $f_{\tau}:\coneH_\tau.y\to [1,\infty)$ as follows
\[
f_{\tau}(z)=\begin{cases} \sum_{w\in I_\tau(z)}\|w\|^{-\alpha} & \text{if $\margI_\tau(z)\neq\emptyset$}\\
\eta_0^{-\alpha}&\text{otherwise}
\end{cases}.
\]

The following proposition supplies an initial dimension which we will bootstrap in the next phase.
Roughly speaking, it asserts that points in $\coneH_{8t}.x_0$ 
(possibly after removing an exponentially small set of exceptions) 
are {\em separated transversal to} $H$, unless $x_0$ is extremely close to a periodic $H$ orbit.    

\begin{proposition}\label{prop:closing lemma intro}
Assume $\Gamma$ is arithmetic. 
There exists $D_1$ (which depends on $\Gamma$ explicitly) satisfying the following. Let $D\geq D_1$ and $x_1\in X$.
Then for all large enough $t$ at least one of the following holds.
 
\begin{enumerate}
\item  There is a subset $I(x_1)\subset [0,1]$ with $|[0,1]\setminus I(x_1)|\ll_X \beta$ 
such that for all $r\in I(x_1)$ we have the following  
\begin{enumerate}
\item $\sfh\mapsto \sfh.a_{8t}u_rx_1$ is injective on $\coneH_{\rws}$.
\item For all $z\in\coneH_\rws.a_{8t}u_rx_1$, we have 
\[
f_{t}(z)\leq  \nuni^{D\rws}.
\]
\end{enumerate}


\item There is $x\in X$ such that $Hx$ is periodic with
\[
\vol(Hx)\leq \nuni^{D_1\rws}\quad\text{and}\quad\dist_X(x,x_1)\leq \nuni^{(-D+D_1)\rws}.
\] 
\end{enumerate} 
\end{proposition}

This proposition is where the arithmeticity assumption on $\Gamma$ is used. The proof is similar to the proof of~\cite[Prop.~6.1]{LM-PolyDensity}.
If we replace the assumption that $\Gamma$ is arithmetic with the weaker requirement that $\Gamma$ has algebraic entries, we get a version of this proposition where part~(2) is replaced with the following. 

\medskip

\begin{enumerate}
    \item[(2')]  {\em There is $x\in X$ with
    \[
\dist_X(x,x_1)\leq \nuni^{(-D+D_1)\rws},
\] 
satisfying the following: there are elements $\gamma_1$ and $\gamma_2$ in ${\rm Stab}_H(x)$
with $\|\gamma_i\|\leq\nuni^{D_1t}$ for $i=1,2$ so that the group generated by $\{\gamma_1,\gamma_2\}$ is Zariski dense in $H$.}
\end{enumerate}

\section{Improving the dimension}
Proposition~\ref{prop:closing lemma intro} shows that up to an exponentially small error, the set $\{a_{9t}u_rx_0: r\in[0,1]\}$ has a {\em small} positive dimension {\em transversal} to $H$ at controlled scales. The objective in this step is to show that by applying elements of the form $a_\ell u_r$, for a fixed $\ell$ and a random $r\in[0,1]$, we can inductively improve this dimension {\em transversal} to $H$ at controlled scales to $\alpha$ (which will be chosen to be smaller than but quite close to $1$). This is achieved by investigating the evolution of a certain Margulis function (cf.\ the survey \cite{EskinMozes-MF}, though the Margulis function we use here is somewhat intricate).

To state the main result which is Proposition~\ref{propos: main bootstrap equi}, we need some notation. Let $R$ be as in Theorem~\ref{thm:main}, and set $t=\frac{1}{D_1}\log R$ with $D_1$ as in Proposition~\ref{prop:closing lemma intro}. We will assume $R$ is large enough so that the conclusion of Proposition~\ref{prop:closing lemma intro} holds with this $t$.

Let 
\[
\coneH=\boxHs_\beta\cdot\Bigl\{u_r: |r|\leq \eta_0\Bigr\}.
\]
where $\beta=\nuni^{-\kappa t}$ for some small $\kappa>0$. (More explicitly, we will fix some $0<\vare\leq10^{-8}$ which will depend on $\mixexp$, and let  $\kappa=10^{-6}\vare/D$ where $D= D_1(D_0+1)$.) 

It will be more convenient to {\em approximate} translations 
$\{a_su_rx_0: r\in[0,1]\}$ with sets which are a disjoint union of local $\coneH$-orbits as we now define.  
Let $F\subset B_\rfrak(0,\beta)$ be a finite set with $\#F\geq \nuni^{t/2}$, and let $y\in X$. For every $w\in F$, let $\coneH_w\subset \coneH$ be a Borel set so that $m_H(\coneH_w)\geq \beta^{4}$ and 
$m_H\Bigl(\coneH_w\triangle(\boxHs_{\beta^{10}}\cdot\coneH_w)\Bigr)\leq \beta m_H(\coneH_w)$,
where $m_H$ denotes a fixed Haar measure on $H$. Put 
\be\label{eq: define cone intro}
\cone=\bigcup\coneH_w.\{\exp(w)y: w\in F\}.
\ee
We equip $\cone$ with the probability measure $\mu_\cone=\frac{1}{\sum_wm_H(\coneH_w)}\sum_{w} \mu_{w,y}$ where 
$\mu_{w,y}$ denotes the pushforward of $m_H|_{\coneH_w}$ to $\coneH_w.\exp(w)y$ for every $w\in F$.

\medskip

Let $\theta$ be a small constant depending on the decay of matrix coefficients in $G/\Gamma$ (the exact value we shall use is $\theta =\mixexp^2/10^4M$, where $\mixexp$ is as~\eqref{eq:actual-mixing-intro} and $M$ as in Proposition~\ref{prop: equi translates of cone intro}). 
Let 
\[
\text{$\balpha=1-\theta\quad$ and $\quad\sqrt\vare=\theta$}.
\]
Let $\ell=0.01\vare t$, and let $\rwm_\ell$ be the probability measure on $H$ defined by 
\[
\rwm_\ell(\varphi)=\ave\varphi(a_{\ell}\uvk)\uvkd\qquad\text{for all $\varphi\in C_c(H)$;}
\] 
let $\rwm_\ell^{(n)}=\rwm_\ell\star\cdots\star\rwm_\ell$ denote the $n$-fold convolution of $\rwm_\ell$ for all $n\in\bbn$.

\medskip

The following is the main statement.

\begin{proposition}\label{propos: main bootstrap equi}
Let $x_1\in X$, and assume that Proposition~\ref{prop:closing lemma intro}(2) does not hold for $D$, $x_1$, and $t$. 
Let $r_1\in I(x_1)$, where $I(x_1)$ is as in Proposition~\ref{prop:closing lemma intro}(1), and put $x_2=\mathsf a_{8t}u_{r_1}x_1$. Let  
\[
J:=[d_{1}-{10^4}\vare^{-1/2}, d_{1}]\cap\bbn,
\] 
where $d_{1}=100\lceil\tfrac{4D-3}{2\vare}\rceil$. 
For every $d\in J$, there is a collection 
\[
\Xi_d=\{\cone_{d,i}: 1\leq i\leq N_d\}
\]
of sets defined as in~\eqref{eq: define cone intro} for some $F_{d,i}\subset B_\rfrak(0,\beta)$, 
so that both of the following hold: 
\begin{enumerate}
    \item Put $\scmf=\nuni^{-\sqrt\vare t}$. Let $d\in J$, $1\leq i\leq N_d$, and let $w_0\in B_\rfrak(0,\beta)$. 
    Then for every $w\in B_\rfrak(w_0,\varrho)$ and all $b\geq \nuni^{-t/2}$, we have
\be\label{eq: energy estimate final intro-1}
\frac{\#\Bigl(B_\rfrak(w,b) \cap B_\rfrak(w_0, \varrho)\cap F_{d,i}\Bigr)}{\#\Bigl(B(w_0, \varrho)\cap F_{d,i}\Bigr)}\leq \nuni^{\vare t} (b/\varrho)^{\alpha}. 
\ee
\item For all $s\leq t$ and all $r\in[0,2]$, we have
\be\label{eq: nud1 and nud1-d}
\biggl|\int \varphi(a_su_rz)\diff\!\rwm_\ell^{(d_{1})}\star\mu_{\coneH_{t}.x_2}(z)-\sum_{d,i}c_{d,i}\int \varphi(a_su_rz)\diff\!\rwm_\ell^{(d_{1}-d)}\star\mu_{\cone_{d,i}}(z)\biggr|
\ll \operatorname{Lip}(\varphi)\beta^{\ref{k: bootstrap beta exp}}
\ee
where $\varphi\in C_c^{\infty}(X)$, $c_{d,i}\geq0$ and $\sum_{d,i}c_{d,i}=1-O(\beta^{\ref{k: bootstrap beta exp}})$, $\operatorname{Lip}(\varphi)$ is the Lipschitz norm of $\varphi$, and $\constk\label{k: bootstrap beta exp}$ and the implied constants depend on $X$.
\end{enumerate}
\end{proposition}

Roughly speaking, the proposition states 
that up to an exponentially small error, 
$\rwm_\ell^{(d_1)}\star\mu_{\coneH_t.x_1}$ may be decomposed as 
$\sum_{d,i}c_{d,i}\rwm_\ell^{(d_1-d)}\star\mu_{\cone_{d,i}}$ where $\sum_{d,i} c_{d,i}=1-O(\beta^{\ref{k: bootstrap beta exp}})$ (see~\eqref{eq: nud1 and nud1-d}) and for all $d\in J$ and $1\leq i\leq N_d$ the dimension of $\cone_{d,i}$ transversal to $H$ at controlled scales is $\geq \alpha$ (see~\eqref{eq: energy estimate final intro-1}).

\medskip

Combining Proposition~\ref{propos: main bootstrap equi} with the previous discussion, we may complete the proof of Theorem~\ref{thm:main}.
A brief outline of this deduction follows: Let $x_0$ be as in the statement and suppose that part~(2) in Theorem~\ref{thm:main} does not hold. 
We first apply Proposition~\ref{prop: linearization translates} with $s=\log T-C\log R$ for an appropriate large constant $C$ to improve the (weak) Diophantine property of $x_0$ provided by the failure of part~(2) in Theorem~\ref{thm:main} to the stronger Diophantine property, 
\begin{equation}\label{eq: stronger diophantine equation}
    d_X(x,x_1)\gg R^{-D_0} \qquad\text{for all $x$ with $\vol(Hx)\leq R$}
\end{equation}
for most points $x_1$ on $\{a_{\log T-C\log R}u_rx_0: r\in[0,1]\}$. 
Thus, in order to show that $\int_0^1 \varphi(a_{\log T}u_rx_0)\diff\!r$ is within $R^{-\star}$ of $\int \phi$ it is enough to show the same for $\int_0^1 \varphi(a_{C\log R}u_rx_1)\diff\!r$ for $x_1$ satisfying the stronger Diophantine property \eqref{eq: stronger diophantine equation}.

The remaining time, i.e., $C\log R$, will be divided into {\em three phases}: 

\subsubsection*{Phase I}
Recall that $t=\frac{1}{D_1}\log R$. We apply Proposition~\ref{prop:closing lemma intro} with the point $x_1$. Then for every $r_1\in I(x_1)$, the conclusion of part~(1) in that proposition holds for $x_2=a_{8t}u_{r_1}x_1$.
That is,
$h\mapsto hx_2$ is injective over $\coneH_t$ and the transverse dimension of $\coneH_t.x_2$ is $\geq 1/D$ for all 
\begin{equation}\label{eq:x2}
    x_2 \in\Bigl\{a_{8t} u_{r_1}  x_1 : r_1 \in I(x_1)\Bigr\}
\end{equation}
where $D=D_0D_1+D_1$.
Therefore, in order to show that $\int_0^1 \varphi(a_{C\log R}u_rx_1)\diff\!r$ is within $R^{-\star}$ of $\int \phi$, it is enough to show a similar estimate for
$\int_0^1 \varphi(a_{C\log R-8t}u_rx_2)\diff\!r$ for all $x_2$ as in \eqref{eq:x2}.

\subsubsection*{Phase II}
Let $s=2\sqrt\vare t$ (note that this is much larger than $\ell=0.01\vare t$). Then
\[
\int_0^1\varphi(a_{s+ d_1\ell+t}u_rx_2)\diff\!r
\]
is within $R^{-\star}$ of  
\[
\int_0^1\int \varphi(a_{s}u_rz)\diff\!\rwm_\ell^{(d_1)}\star\mu_{\coneH_{t}.x_2}(z)\diff\!r.
\]

We now use Proposition~\ref{propos: main bootstrap equi} to improve the small transversal dimension from $1/D$ to $\alpha$.
More precisely, Proposition~\ref{propos: main bootstrap equi} shows that 
\[
\int_0^1\int \varphi(a_{s}u_rz)\diff\!\rwm_\ell^{(d_{1})}\star\mu_{\coneH_{t}.x_2}(z)\diff\!r
\]
is within $R^{-\star}$ of a convex combination of integrals of the form 
\[
\int_0^1\int \varphi(a_su_rz)\diff\!\rwm_\ell^{(n)}\star\mu_{\cone}(z)\diff\!r
\]
where $0\leq n=d_1-d\leq 10^4\vare^{-1/2}$ and $\cone=\cone_{d,i}$ has dimension at least $\alpha$ transversal to $H$ at controlled scales, see~\eqref{eq: energy estimate final intro-1}. 

\subsubsection*{Phase III}
It now suffices to show that $\int_0^1\int \varphi(a_su_rz)\diff\!\rwm_\ell^{(n)}\star\mu_{\cone}(z)\diff\!r$ is within $R^{-\star}$ of $\int\varphi$ for all $\cone$ and $n$ as above. We will use Proposition~\ref{prop: equi translates of cone intro} to show this. 
First note that   
\[
\int_0^1\int \varphi(a_su_rz)\diff\!\rwm_\ell^{(n)}\star\mu_{\cone}(z)\diff\!r
\]
is within $R^{-\star}$ of 
\[
\int\int_0^1\varphi(a_{s+n\ell}u_rz)\diff\!r\diff\mu_{\cone}(z).
\]
Moreover, we have 
\[
2\sqrt\vare t\leq s+n\ell\leq 2\sqrt\vare t+\frac{10^4\ell}{\sqrt\vare}=102\sqrt\vare t;
\]
in view of our choice of $\theta$ the right most term in the above series of inequalities is $\leq \frac{\mixexp^2}{4M\theta}\sqrt\vare t$.
Thus, Proposition~\ref{prop: equi translates of cone intro}, applied with $\theta=1-\alpha$, $c=2\vare$, $\varrho=\nuni^{-\sqrt\vare t}$, $\rho=\nuni^{-t/2}$, and $\tau=s+n\ell$, gives
\be\label{eq: using prop equi trans}
\biggl|\iint \varphi(a_{s+n\ell}u_rz)\diff\!\mu_{\cone}(z)\diff\!r-\int\varphi \diff\!m_X\biggr|\ll\Sob(\varphi)e^{-\star t}=\Sob(\varphi)R^{-\star}
\ee
where the implied constants depend on $X$. 

Note that the total time required for these three phases is $s+d_1\ell+9t$ which in view of the choices of $s$, $\ell$ and $t$ is indeed a (large) constant times $\log R$. Theorem~\ref{thm:main} follows. 

\medskip

In the setting of Theorem~\ref{thm: main arithmeticity relaxed}, we cannot utilize Proposition~\ref{prop: linearization translates} combined with Proposition~\ref{prop:closing lemma intro} as we did above; cf.\ the weaker conclusion in (2') following Proposition~\ref{prop:closing lemma intro}. Thus, we only use the three phases above (with $t=\star\delta\log T$) to improve the small dimension, namely the parameter $\delta$ in Theorem~\ref{thm: main arithmeticity relaxed}, to $\alpha$. Thus the number of steps required is $\gg 1/\delta$ which forces $\kappa\ll \delta$. Hence, we only obtain the rate $T^{-\star\delta^2}$ in part~(1) of Theorem~\ref{thm: main arithmeticity relaxed}.

\medskip

The proof of Propositions~\ref{propos: main bootstrap equi} relies on the evolution of the Margulis function $\mfm_{\cone,\scmf,\trct}$ defined below. For every $(h,z)\in H\times \cone$, put 
\[
\margI_{\cone,\scmf}(h,z):=\Bigl\{w\in \rfrak: 0<\|w\|<\scmf\text{ and } \exp(w) h z\in h\cone\Bigr\}.
\]
Since $\mathsf E$ is bounded, $\margI_{\cone,\scmf}(h,z)$ is a finite set for all $(h,z)\in H\times\cone$. 

For every $\trct\geq 0$, define the modified 
and localized Margulis function 
$
\mfm_{\cone,\scmf,\trct}: H\times \cone\to [1,\infty)
$
as follows\footnote{If $\Gamma$ is nonuniform, this Margulis function needs to be modified accordingly; but in this section we limit ourselves to the compact case.}: if $\#\margI_{\cone,\scmf}(h,z)\leq \trct$, put 
\[
\mfm_{\cone,\scmf,\trct}(h,z)=\scmf^{-\alpha},
\]
and if $\#\margI_{\cone,\scmf}(h,z)> \trct$, put
\[
\mfm_{\cone,\scmf,\trct}(h,z)=\min\left\{\sum_{w\in I}\|w\|^{-\alpha}: \begin{array}{c}I\subset I_{\cone,\scmf}(h,z)\text{ and }\\ \#(I_{\cone,\scmf}(h,z)\setminus I)=\trct\end{array}\right\}.
\]

\medskip

We begin the outline of the proof of Propositions~\ref{propos: main bootstrap equi} with the following observation: the set $\coneH_t.x_1$ gives rise to sets $\cone$ which are defined as in~\eqref{eq: define cone intro} and since we assume that Proposition~\ref{prop:closing lemma intro}(2) does not hold these satisfy 
$\mfm_{\cone, \scmf, 0}(e,z)\leq \nuni^{Dt}$
for all $z\in\cone$. 

Let $\cone$ be one of these sets, then up to an exponentially small error, $\rwm_\ell\star\mu_\cone$ may be decomposed into $\sum c_j'\mu_{\cone'_j}$ where $c'_j\geq 0$ and $\sum c'_j=1-O(\beta^\star)$. This can be seen by first decomposing $\rwm_\ell\star\mu_\cone$ into a combination of measures supported on subsets of $X$ which are exponentially thin in the direction of $U^-$ (note that $a_\ell\coneH.y$ will be exponentially thin in the direction of $U^-$), 
and then smearing these measures with $B^H_\beta$. Continuing inductively, we obtain the decomposition~\eqref{eq: nud1 and nud1-d}.

The fact that the energy estimate~\eqref{eq: energy estimate final intro-1} is also satisfied for the terms appearing in~\eqref{eq: nud1 and nud1-d} is at the heart of the proposition. The proof of this fact is based on the following bootstrap step: for all but an exponentially small set of $r\in[0,1]$ and all but an exponentially small set of $z\in\cone$, 
\be\label{eq: one step of walk}
\mfm_{\cone,\scmf, \beta^{-1}}(a_{\ell}u_r,z)\leq \nuni^{-\ell/2}\nuni^{Dt}+4\nuni^{\ell}\scmf^{-\alpha}\noI_{\cone, \scmf}(a_{\ell}u_r,z).
\ee 
where $\noI_{\cone,\scmf}(h,z):=\max\bigl\{\#\margI_{\cone,\scmf}(h,z),1\bigr\}$.

The proof of~\eqref{eq: one step of walk} relies on the aforementioned works of Wolff, Schlag, and Zahl~\cite{Wolff, Schlag, Zahl}, see also the discussion below.  

The set $a_\ell u_r.\cone$ can now be used to construct sets $\cone'$ (defined as in~\eqref{eq: define cone intro}), and in view of~\eqref{eq: one step of walk}, the estimate for $\mfm_{\cone',\scmf, \beta^{-1}}$ is {\em improved}. Continuing inductively and using 
$Dt-0.5d_1\ell\leq 3t/4$, after $d\leq d_1$ steps, we have  
\be\label{eq: energy estimate final intro-2}
\mfm_{\cone',\scmf, d\beta^{-1}}(e,z)\leq \nuni^{\vare t}\scmf^{-\alpha}\noI_{\cone',\scmf}(e,z),
\ee
which implies the dimension estimate~\eqref{eq: energy estimate final intro-1} for the set $F'\subset B_\rfrak(0,\beta)$ which is used in the definition of $\cone'$, see~\eqref{eq: define cone intro}.  

\medskip

We emphasize that our inductive scheme produces sets $\cone'$ at every step $0<d\leq d_1$ with an improved bound on $\mfm_{\cone',\varrho, d\beta^{-1}}$, however, it does not guarantee that~\eqref{eq: energy estimate final intro-2} is only satisfied for $d\in J$. On the other hand, and as it was discussed above, the fact that our {\em stopping times} $d$ satisfy $d_1-d\leq 10^4\vare^{-1/2}$ is essential for us when we apply Proposition~\ref{prop: equi translates of cone intro}. We remedy this issue as follows: if~\eqref{eq: energy estimate final intro-1} holds for some $\cone'$ defined at step $d<d_1-10^4\vare^{-1/2}$, then we use the above inductive scheme to show that starting with $\cone'$, the process {\em again} terminates (i.e.,~\eqref{eq: energy estimate final intro-2} is satisfied) in at most $10^4\vare^{-1/2}$-many further steps.

On a related note, it should also be mentioned that ideally one would like to replace the interval of possible choices of $d\in J$ with the singleton $\{d_{1}\}$, i.e., to show that after exactly $d_{1}$ steps, one obtains sets which satisfy~\eqref{eq: energy estimate final intro-1}. Indeed, such statement can be obtained if one is content with restricting to $\alpha<1/2$ --- this can be achieved using estimates analogues to~\cite[Lemma 5.1]{EMM-Upp} where one replaces the integral over $[0,1]$ with integrals over {\em much smaller} intervals and by conditioning the random walk. 

However, it is essential for us to work with $\alpha=1-\theta$ where $\theta>0$ is rather small. 
For such choices of $\alpha$, there are vectors $w\in\rfrak$ where the growth of $\|a_tu_rw\|$ is too slow. 
Indeed, in general, one can only guarantee that 
\be\label{eq: linear algebra intro}
\ave\|a_{t}{\uvk} w\|^{-\alpha}\uvkd\ll\nuni^{-\theta t}{\|w\|^{-\alpha}}.
\ee
Using this general fact (which is an exercise in linear algebra) as an input, one can prove a version of Proposition~\ref{propos: main bootstrap equi} where 
\[
\text{$[d_{1}-\tfrac{10^4}{\sqrt\vare}, d_{1}]\quad$ is replaced by $\quad[d_1-\tfrac{36}{\theta\sqrt\vare}, d_1]$}, 
\]
in particular, the length of the interval cannot be made smaller than $(\theta\sqrt\vare)^{-1}$.

As it was discussed above, this improvement is pivotal to our analysis,
and it is made possible by bringing to bear~\cite{Wolff, Schlag, Zahl} in this step of the argument as well. Indeed, the {\em poor} rate in~\eqref{eq: linear algebra intro} is closely related to the existence of double zeroes for the map $r\mapsto (a_tu_rw)_{12}$ (the $(1,2)$-entry of $a_tu_rw$); we thus use~\cite{Wolff, Schlag, Zahl}, to control {\rm tangencies} 
between the two parabolas $\{(u_rw_1)_{12}: r\in[0,1]\}$ and $\{(u_rw_2)_{12}: r\in[0,1]\}$ for most pairs $w_1,w_2\in F$ and most $r\in [0,1]$. 
This yields an improved version of~\eqref{eq: linear algebra intro} which we use crucially.

\bibliographystyle{amsalpha}
\bibliography{papers}

\end{document}